\documentclass[11pt,a4wide,oneside,reqno]{amsart}
\usepackage{amssymb,amsmath,url,a4wide}
\usepackage[utf8]{inputenc}
\usepackage[left=3cm,right=3cm,top=3cm,bottom=2cm]{geometry}
\usepackage{enumitem}

\newcommand{\citep}[1]{\cite{#1}}

\newcommand{\F}{\mathbb{F}}

\newcommand{\R}{\mathbb{R}}
\newcommand{\E}{\mathbb{E}}
\newcommand{\A}{\mathcal{A}}
\newcommand{\e}{\varepsilon}

\newcommand{\rv}{\vert}
\newcommand{\lv}{\vert}

\setlist[enumerate]{label=(\roman*)}

\usepackage[breaklinks=true]{hyperref}
\usepackage{cleveref}

\numberwithin{equation}{section}
\newtheorem{theorem}{Theorem}[section]

\newtheorem{corollary}[theorem]{Corollary}

\theoremstyle{definition}

\usepackage{thmtools}
\declaretheoremstyle[%
  spaceabove=0pt,%
  spacebelow=0pt,%
  headfont=\normalfont\itshape,%
  postheadspace=1em,%
  qed=\qedsymbol%
]{mystyle}

\usepackage[url=false,doi=true,maxnames=4,giveninits=true,isbn=false]{biblatex}
\usepackage{xpatch}
\xpatchbibmacro{volume+number+eid}{%
  \setunit*{\adddot}%
}{%
}{}{}
\DeclareFieldFormat[article]{number}{\mkbibparens{#1}}
\renewbibmacro{in:}{%
  \ifentrytype{article}{}{\printtext{\bibstring{in}\intitlepunct}}}
\bibliography{biblio}

\title{Strongly common graphs with odd girth are cycles}

\author{Leo Versteegen}
\thanks{The author is grateful to be funded by Trinity College of the University of Cambridge through the Trinity External Researcher Studentship.}
\address{Department of Pure Mathematics and Mathematical Statistics, Centre for Mathematical Sciences, Wilberforce Road, Cambridge CB3 0WB, United Kingdom}
\email{lvv23@dpmms.cam.ac.uk}

\begin{document}
\maketitle
\begin{abstract}
A graph $H$ is called \emph{strongly common} if for every coloring $\phi$ of $K_n$ with two colors, the number of monochromatic copies of $H$ is at least the number of monochromatic copies of $H$ in a random coloring of $K_n$ with the same density of color classes as $\phi$. In this note we prove that if a graph has odd girth but is not a cycle, then it is not strongly common. This answers a question of Chen and Ma.
\end{abstract}
\section{Introduction}
Let $H=(V(H),E(H))$ be a graph with $e(H)$ edges on $v(H)$ edges. For a probability space $(X,\A,\mu)$ and an integrable function $\phi\colon X^2\rightarrow \R$, we define the \emph{homomorphism density of $H$ in $\phi$} as
\begin{align*}
t_H(\phi)=\E_{x\in X^V(H)} \prod_{vw\in E(H)} \phi(x_v,x_w),
\end{align*}
where the expectation is taken with respect to the $v(H)$-fold product measure $\mu^{\otimes v(H)}$. The graph $H$ is called \emph{common}\footnote{The original definition of ``common'' counts the monochromatic copies of $H$ in a coloring of $K_n$ as $n$ goes to infinity, but for us the definition via graphons will be more convenient.} if for every \emph{graphon}, i.e., for every measurable map $W\colon [0,1]^2\rightarrow [0,1]$, the inequality $t_H(W)+t_H(1-W)\geq 2^{1-e(H)}$ holds. Graphs with this property have been studied extensively since the late 1950s when Goodman first pointed out that triangles possess it \cite{goodman59}, but it was only in 2022 that Behague, Morrison and Noel introduced the following related notion \cite{behague2022common}. A graph is called \emph{strongly common}, if 
\begin{align*}
t_H(W)+t_H(1-W)\geq t_{K_2}(W)^{e(H)}+t_{K_2}(1-W)^{e(H)}
\end{align*}
for every graphon $W$. Beyond those graphs that are trivially strongly common because they have Sidorenko's property, Behague, Morrison and Noel showed that the triangle and the cycle of length five are strongly common. Furthermore, they conjectured that any odd cycle is strongly common, which has since been confirmed by Kim and Lee \cite{kim2022cycles}.

In the other direction, Chen and Ma showed that any graph $H$ containing a triangle that is not itself a triangle is not strongly common \cite{chen2023triangle}. They also asked whether their result could be generalized in the sense that if a graph has odd girth $k$ but is not isomorphic to $C_k$, then it cannot be strongly common. In the present note, we give an affirmative answer to their question, and remarkably, the same, extremely simple function may be used as a witness for all $H$. In what follows, we take $\mu$ to be the uniform probability measure on $\F_2$.

\begin{theorem}\label{thm:main}
There exists a function $f\colon \F_2^2\rightarrow \{\pm 1\}$ with $\E(f)=0$ such that for any graph $H$ that has odd girth and is not a cycle, and every $\alpha\in (0,1)\setminus \{1/2\}$, there exists $\e_0>0$ such that for all $\e<\e_0$, we have
\begin{align*}
t_H(\alpha+\e f)+t_H(1-\alpha-\e f)<\alpha^{e(H)}+(1-\alpha)^{e(H)}.
\end{align*} 
\end{theorem}
Using any measure-preserving map $[0,1]\rightarrow \F_2$, the function $\alpha+\e f$ in the theorem above may be converted into a graphon, which yields the following corollary.
\begin{corollary}
If $H$ is a graph that has odd girth and is not a cycle, then $H$ is not strongly common.
\end{corollary}
We remark that together with the result of Kim and Lee that odd cycles are strongly common \cite{kim2022cycles}, this completely characterizes the strongly common graphs with odd girth. For graphs with even girth on the other hand, the problem is still completely open. In particular, it would be interesting to find a strongly common bipartite graph for which it is not known whether it has Sidorenko's property.
\section{Proof of the theorem}
Our proof of \Cref{thm:main} requires only one additional piece of notation. For a graph $H$ containing a vertex $u$ and a vector $x\in \F_2^{V(H)}$, we denote by $\hat{x}_u$ the vector in $\F_2^{V(H)\setminus\{u\}}$ that is obtained from $x$ by omitting the coordinate $x_u$.
\begin{proof}[Proof of \Cref{thm:main}]
The function $f$ is given by $f(x,y)=(-1)^{x+y+1}$. Let $H$ and $\alpha$ be as in the claim and denote by $k$ the girth of $H$, which is odd by assumption. We may assume without loss of generality that $\alpha >1/2$ and take $\e_0=\min(1-\alpha, 2^{-2e(H)}(\alpha-1/2))$. Writing $\phi(x,y)$ for $\alpha+\e f(x,y)$, we have
\begin{align}\label{eq:overall-bound}
t_H(\phi)&=\E_{x\in \F_2^{V(H)}} \prod_{vw\in E(H)} (\alpha+\e f(x_v,x_w))\nonumber \\
&= \sum_{S\subseteq E(H)} \e^{\lv S\rv}\alpha^{e(H)-\lv S\rv} \E_{x\in \F_2^{V(H)}} \prod_{vw\in S} f(x_v,x_w).\nonumber\\
&\leq \sum_{\substack{S\subset E(H)\\ \lv S\rv \leq k}}\left[ \e^{\lv S\rv}\alpha^{e(H)-\lv S\rv} \E_{x\in \F_2^{V(H)}} \prod_{vw\in S} f(x_v,x_w)\right] + 2^{e(H)}\e^{k+1} .
\end{align}
Suppose now that the edges in a given non-empty set $S\subset E(H)$ with $\lv S\rv\leq k$ do not form a cycle. Since $H$ has girth $k$, there must be a vertex $u\in V(H)$ that is incident to an odd number of edges in $S$. Thus, there must exist a map $r\colon \F_2^{V(H)\setminus \{u\}}\rightarrow \{\pm 1\}$ such that for all $x\in \F_2^{V(H)}$,
\begin{align*}
\prod_{vw\in S} f(x_v,x_w)=(-1)^{x_u+r(\hat{x}_u)}.
\end{align*}
From this it follows that
\begin{align*}
\E_{x\in \F_2^{V(H)}} \prod_{uv\in S} f(x_u,x_v)=\E_{x_u\in \F_2} (-1)^{x_u}\E_{x'\in \F_2^{V(H)\setminus\{u\}}} (-1)^{r(x')}=0.
\end{align*}
On the other hand, if the edges of $S$ do form a cycle (implying that $\lv S\rv=k$), it is easy to see that
\begin{align*}
\E_{x\in \F_2^{V(H)}} \prod_{uv\in S} f(x_u,x_v)=(-1)^{\lv S\rv}=-1.
\end{align*}
Inserting this information into \eqref{eq:overall-bound}, we obtain
\begin{align*}
t_H(\phi)\leq \alpha^{e(H)}-\e^k\alpha^{e(H)-k} \lv \{S\subset E(H):S\cong C_k\}\rv + 2^{e(H)}\e^{k+1}.
\end{align*}
By a similar calculation, we see that
\begin{align*}
t_H(1-\phi)\leq (1-\alpha)^{e(H)}+\e^k(1-\alpha)^{e(H)-k} \lv \{S\subset E(H):S\cong C_k\}\rv + 2^{e(H)}\e^{k+1}.
\end{align*}
Combining these two bounds and the fact that $H$ has at least one cycle of length $k$ yields
\begin{align*}
t_H(\phi)+t_H(1-\phi)\leq \alpha^{e(H)}+(1-\alpha)^{e(H)} +\e^k (-\alpha^{e(H)-k} +(1-\alpha)^{e(H)-k})+ 2^{e(H)+1}\e^{k+1}.
\end{align*}
Using the mean value theorem, we deduce from the above that
\begin{align*}
t_H(\phi)+t_H(1-\phi)\leq \alpha^{e(H)}+(1-\alpha)^{e(H)} -2^{k-e(H)}\e^k (\alpha-1/2)+ 2^{e(H)+1}\e^{k+1},
\end{align*}
which is less than $\alpha^{e(H)}+(1-\alpha)^{e(H)}$ as $\e<2^{-2e(H)}(\alpha-1/2)$, completing the proof.
\end{proof}

\printbibliography

\end{document}